\documentclass[12pt,reqno]{amsart}
\usepackage{comment}
\usepackage{amsmath,amsfonts,amsthm}
\usepackage{tikz}
\newtheorem{theorem}{Theorem}[section]
\newtheorem{prop}[theorem]{Proposition}
\newtheorem{lemma}[theorem]{Lemma}

\newtheorem{definition}[theorem]{Definition}

\usepackage[a4paper,
top=24mm,bottom=24mm,left=26mm,right=26mm
]{geometry}

\usepackage[numbers,sort]{natbib}
\usepackage[
bookmarks=true,
bookmarksdepth=subsection,
bookmarksnumbered=true,
breaklinks=true,
colorlinks=true,
linkcolor=blue,
citecolor=magenta
]{hyperref}

\usepackage[]{graphicx}
\usepackage{booktabs}
\usepackage{tikz}
\usetikzlibrary{cd}
\usetikzlibrary{arrows.meta}
\usetikzlibrary{matrix,arrows,decorations.pathmorphing}
\tikzcdset{arrow style=tikz,
  diagrams={>=Straight Barb}
}
\usetikzlibrary{positioning}
\usetikzlibrary{decorations.markings}
\tikzset{->-/.style={decoration={
      markings,
      mark=at position #1 with {\arrow{>}}}, postaction={decorate}}}
\tikzset{->>-/.style={decoration={
      markings,
      mark=at position #1 with {\arrow{>>}}}, postaction={decorate}}}

\usepackage{mathtools}
\usepackage{makecell}
  
\usepackage{mathrsfs}
\usepackage[sb]{libertine}
\usepackage[varqu,varl,scale=0.9]{zi4}
\usepackage[libertine,bigdelims,vvarbb]{newtxmath}
\numberwithin{equation}{section}
\newcommand{\I}{\mathrm{i}}

\DeclareMathOperator{\sh}{sh}
\DeclareMathOperator{\ch}{ch}
\DeclareMathOperator{\End}{End}
\DeclareMathOperator{\Mod}{Mod}
\DeclareMathOperator{\Sk}{Sk}

\DeclareMathDelimiter{\Norm}{\mathord}{largesymbols}{"3E}{largesymbols}{"3E}

\allowdisplaybreaks
\begin{document}

\baselineskip 15pt
\parskip 7pt
\sloppy


\title[]{Generalized Double Affine Hecke Algebra for Double Torus}
\dedicatory{To the memory of my father}


\author[
]
{Kazuhiro Hikami}

\address{Faculty of Mathematics,
  Kyushu University,
  Fukuoka 819-0395, Japan.}

\email{
  \texttt{khikami@gmail.com}
}



\date{February 14, 2024}

\begin{abstract}
  We propose a generalization of the double affine Hecke algebra of
  type-$C^\vee C_1$ at specific parameters by introducing a ``Heegaard
  dual'' of the Hecke operators.
  Shown is  a relationship with the skein algebra on double torus.
  We give automorphisms of the algebra
  associated with the Dehn twists on
  the double torus.
\end{abstract}


\keywords{double affine Hecke algebra, skein algebra,
  Askey--Wilson polynomial}

\subjclass[2000]{
}


\maketitle
\section{Introduction}

The double affine Hecke algebra
(DAHA)
was introduced by Cherednik for
studies on the Knizhnik--Zamolodchikov equation and its applications
to  orthogonal symmetric polynomials via the Dunkl--Cherednik
operators~\cite{Chered05Book}.
It is a fundamental modern tool in
mathematics and  physics.
One of  applications of DAHA
is the skein algebra,
which receives renewed interests from a viewpoint of the quantization
of the character varieties.
Although, the DAHA  is so far only applicable for the skein
algebra on the
once-punctured torus~$\Sigma_{1,1}$
and the 4-punctured sphere~$\Sigma_{0,4}$.
The former is the DAHA of type-$A_1$,
and the latter is of
type-$C^\vee C_1$~\cite{Oblom04a,BerestSamuel16a,KHikami19a,MortoSamue21a}.
The superpolynomials for links  on $\Sigma_{1,1}$ and $\Sigma_{0,4}$
were constructed by use of 
the automorphism of DAHA~\cite{IChered13a,IChered16a}.

A generalization of DAHA
for the skein algebra on  higher-genus surface  was initiated
in~\cite{ArthaShaki15a,ArthaShaki17a}.
Constructed are $q$-difference operators for generators of the skein
algebra on $\Sigma_{2,0}$ as  a $\mathbb{Z}_3$ generalization
of the DAHA of type-$A_1$.
Though the Iwahori--Hecke algebraic structure seems to be missing,
shown was that
it is isomorphic to the skein
algebra on~$\Sigma_{2,0}$~\cite{CookeSamue21a}.

Another representation of the skein algebra on $\Sigma_{2,0}$
was given in our previous paper~\cite{KHikami19a}.
Combining 
the DAHAs of type-$A_1$ and $C^\vee C_1$,
we constructed the
$q$-difference operators for generators of the skein
algebra~\cite{ArthaShaki15a,ArthaShaki17a}.
Using the automorphisms of $A_1$ DAHA, we
computed the DAHA polynomial for double twist knots, and observed
a relationship with the colored Jones polynomials.
The present paper relies on~\cite{KHikami19a},
but  we rather aim to generalize the $C^\vee C_1$ DAHA
at  specific parameters by
introducing 
``Heegaard dual'' operators of $C^\vee C_1$ DAHA.
Thus our generalization is different from~\cite{EtinObloRain07a}
where
the generalized DAHA was 
associated with a 2-dimensional crystallographic group.

This paper is organized as follows.
In Section~\ref{sec:DAHA_C}, we recall the skein algebra on surface
and
DAHA.
We review an isomorphism between the skein algebra on the 4-punctured
sphere and DAHA of type-$C^\vee C_1$.
In Section~\ref{sec:Sigma20},
we pay attention to $C^\vee C_1$ DAHA at $\mathbf{t}_\star$.
By introducing
Heegaard dual operators,
we propose a generalized DAHA.
Using the automorphisms we study a relationship with the skein
algebra on~$\Sigma_{2,0}$.

Throughout this article
we use
\begin{equation}
  \label{ch_sh}
  \ch(x)=x+x^{-1},
  \qquad \qquad
  \sh(x)=x-x^{-1}.
\end{equation}

\section{Preliminaries on DAHA of \texorpdfstring{$C^\vee C_1$}{C
    C1}-type and Skein Algebra}
\label{sec:DAHA_C}
\subsection{DAHA of \texorpdfstring{$C^\vee C_1$}{C C1}-type}
We recall the
double affine Hecke algebra~${H}_{q,\mathbf{t}}$ of type-$C^\vee C_1$
with 4 parameters
$\mathbf{t}=(t_0,t_1,t_2,t_3)$~\cite{NoumiStokm00a} (see
also~\cite{Chered05Book,Macdonald03book}).
The DAHA $H_{q,\mathbf{t}}$ is
generated by $\mathsf{T}_0^{\pm 1}$, $\mathsf{T}_1^{\pm 1}$,
$(\mathsf{T}^\vee_0)^{\pm 1}$, and
$(\mathsf{T}^\vee_1)^{\pm 1}$ satisfying
\begin{equation}
  \label{eq:17}
  \begin{alignedat}{2}
    &
      \mathsf{T}_0 + t_0- t_0^{-1} =
      \mathsf{T}_0^{-1} ,
      & \qquad \qquad
    &
      \mathsf{T}_1 + t_1- t_1^{-1}=
      \mathsf{T}_1^{-1},
    \\
    &
      \mathsf{T}^\vee_0 + t_2- t_2^{-1}=
      \left(\mathsf{T}^\vee_0 \right)^{-1},
      &&
         \mathsf{T}^\vee_1 +t_3- t_3^{-1} =
         \left( \mathsf{T}^\vee_1 \right)^{-1},
  \end{alignedat}
\end{equation}
and
\begin{equation}
  \label{prod_TTTT}
  \mathsf{T}_1^\vee  \mathsf{T}_1  \mathsf{T}_0  \mathsf{T}_0^\vee
  =q^{-\frac{1}{2}}.
\end{equation}
Here and hereafter we use
\begin{gather}
  \label{Hecke_X}
  \mathsf{X}= \left( \mathsf{T}_1^\vee \mathsf{T}_1 \right)^{-1}
  =
  q^{\frac{1}{2}} \mathsf{T}_0 \mathsf{T}_0^\vee
  ,
  \\
  \label{Hecke_Y}
  \mathsf{Y}=\mathsf{T}_1  \mathsf{T}_0 .
\end{gather}
The spherical  DAHA
is defined by
\begin{equation}
  \label{spherical-DAHA}
  SH_{q,\mathbf{t}}
  =
  \mathsf{e} \, H_{q,\mathbf{t}} \, \mathsf{e} , 
\end{equation}
where the idempotent $\mathsf{e}$ is
\begin{equation}
  \label{idempotent_AW}
  \mathsf{e} = \frac{1}{t_1+t_1^{-1}}\left( t_1 + \mathsf{T}_1 \right)
  .
\end{equation}

The polynomial representation is given as~\cite{NoumiStokm00a}
\begin{equation}
  \label{AW_poly_rep}
  \begin{aligned}
    \mathsf{T}_0
    & \mapsto
      {t_0}^{-1} \, \mathsf{s} \, \eth -
      \frac{
      q^{-1} \left({t_0}^{-1}-t_0\right)  x^2
      + q^{-\frac{1}{2}} \left({t_2}^{-1}-t_2 \right)  x
      }{
      1-q^{-1} x^2
      } \,
      \left( 1 - \mathsf{s} \, \eth \right)
      ,
    \\
    \mathsf{T}_1 
    & \mapsto
      {t_1}^{-1} \mathsf{s}+
      \frac{
      \left( {t_1}^{-1}-t_1 \right)
      +
      \left( {t_3}^{-1}-t_3 \right)  x
      }{
      x^2-1} \,
      \left( \mathsf{s} -1 \right)
      ,
    \\
    \mathsf{T}_0^\vee
    & \mapsto
      q^{-\frac{1}{2}} {\mathsf{T}_0}^{-1} x
      ,
    \\
    \mathsf{T}_1^\vee
    & \mapsto
      x^{-1} {\mathsf{T}_1}^{-1}
      .
  \end{aligned}
\end{equation}
Here we mean
\begin{equation}
  \left(\mathsf{s} \, f\right)(x)= f(x^{-1}),
  \qquad\qquad\qquad
  \left(\eth \, f\right)(x)= f(q \,  x) .
\end{equation}
Then  the idempotent~\eqref{idempotent_AW} is a map to the symmetric
Laurent polynomials,
$ \mathsf{e}:  \mathbb{C}[x^{\pm 1}]  \to  \mathbb{C}[x+x^{-1}] $.
The Askey--Wilson operator is given from~\eqref{Hecke_Y} as
\begin{equation}
  \label{difference_Askey-Wilson}
  \left.
    \ch(\mathsf{Y} )
  \right|_{\text{sym}}
  \mapsto
  W(x;\mathbf{t}) \left( \eth-1 \right)
  +W(x^{-1}; \mathbf{t}) \left(\eth^{-1} -1\right)
  + t_0 t_1+ \left(t_0 \,  t_1\right)^{-1}
  ,
\end{equation}
where $\text{sym}$ denotes an action on the symmetric Laurent
polynomials, and
\begin{equation}
  \label{eq:14}
  W(x;\mathbf{t})
  =
  t_0 \, t_1 \,
  \frac{\left( 1- \frac{1}{t_1t_3} x \right)
    \left( 1+\frac{t_3}{t_1}x \right)
    \left( 1- q^{\frac{1}{2}} \frac{1}{t_0t_2} x\right)
    \left( 1+ q^{\frac{1}{2}} \frac{t_2}{t_0} x\right)
  }{
    \left(1-x^2\right)
    \left( 1- q  \, x^2 \right) 
  } .
\end{equation}
The eigen-polynomial of~\eqref{difference_Askey-Wilson}
is the
Askey--Wilson polynomial $P_m(x;q,\mathbf{t})$,
\begin{equation}
  \label{AW_polynomial}
  P_m(x;q,\mathbf{t})
  =
  \frac{
    (a b, a c, a d; q)_m}{
    a^m (a b c d q^{m-1};q)_m}
  {}_4\phi_3
  \left[
    \begin{matrix}
      q^{-m}, q^{m-1} a b c d, a x, a x^{-1}
      \\
      a b, a c, a d
    \end{matrix}
    ;
    q,q
  \right],
\end{equation}
satisfying
\begin{equation}
  \label{eq:21}
  \left( \mathsf{Y} + \mathsf{Y}^{-1} \right) P_m(x;q,\mathbf{t})
  =
  \ch \left( t_0 t_1 q^{-m}\right) \, P_m(x;q,\mathbf{t}) ,  
\end{equation}
where
\begin{equation*}
  a=\frac{1}{t_1t_3}, 
  \qquad
  b=-\frac{t_3}{t_1} ,
  \qquad
  c=\frac{q^{\frac{1}{2}}}{t_0 t_2},
  \qquad
  d=-\frac{q^{\frac{1}{2}} t_2}{t_0} .
\end{equation*}
See~\cite{AsWi85,GaspRahm04} for properties of the Askey--Wilson polynomials.

\subsection{Skein Algebra on
  \texorpdfstring{$\Sigma_{0,4}$}{Sigma0,4}}

The skein algebra $\Sk_A(\Sigma)$
on surface $\Sigma$ is generated by isotopy classes of framed links in
$\Sigma\times[0,1]$
satisfying
the skein relation
\begin{gather}
  \vcenter{\hbox{
    \begin{tikzpicture}
      \draw [line width=1.2pt](0,1) --( 1,0) ;
      \draw[line width=10pt,white](0,0)--(1,1);
      \draw[line width=1.2pt](0,0)--(1,1);
    \end{tikzpicture}
  }}
  =
  A \,
  \vcenter{\hbox{
    \begin{tikzpicture}
      \draw [line width=1.2pt] (1,1) to[out=-135,in=135] (1,0);
      \draw [line width=1.2pt] (0,1) to[out=-45,in=45] (0,0);
    \end{tikzpicture}
  }}
  +A^{-1} \,
  \vcenter{\hbox{
    \begin{tikzpicture}
      \draw [line width=1.2pt] (1,1) to[out=-135,in=-45] (0,1);
      \draw [line width=1.2pt] (1,0) to[out=135,in=45] (0,0);
    \end{tikzpicture}
  }}
  ,
  \\[2mm]
  \vcenter{\hbox{
    \begin{tikzpicture}
      \draw [line width=1.2pt] (0,0) circle (0.5) ;
    \end{tikzpicture}
  }}
  = -A^2 - A^{-2} .
  \notag
\end{gather}
A multiplication
$\mathbb{x} \, \mathbb{y}$
of links $\mathbb{x}$ and $\mathbb{y}$
means
that $\mathbb{x}$ is vertically above $\mathbb{y}$,
\begin{gather*}
  \mathbb{x} \, \mathbb{y}
  =
  \newcolumntype{C}{>{$}c<{$}}
  \begin{tabular}{|C|}
    \hline
    \hphantom{aa}\mathbb{x}\hphantom{aa}
    \\
    \hline
    \mathbb{y}
    \\
    \hline
  \end{tabular}
\end{gather*}
When two simple closed curves $\mathbb{x}$ and $\mathbb{y}$ on
$\Sigma$
intersect exactly once, we have
\begin{equation}
  \label{example_skein}
  \frac{1}{A^{\pm 2}-A^{\mp 2}}
  \left(A^{\pm 1} \mathbb{x} \, \mathbb{y} - A^{\mp 1} \mathbb{y} \, \mathbb{x}\right)
  =
  \mathscr{D}_{\mathbb{x}}^{\mp 1}(\mathbb{y}) ,
\end{equation}
where
$\mathscr{D}_{\mathbb{x}}$ denotes the left Dehn twist
along~$\mathbb{x}$.
It is noted that
\begin{equation}
  \mathscr{D}_{\mathbb{y}}(\mathbb{x})
  =
  \mathscr{D}_{\mathbb{x}}^{-1}(\mathbb{y}) .
\end{equation}
A finite
set of Dehn twists along non-separating simple closed curves
generates
the mapping class group $\Mod(\Sigma_{g,0})$
of a surface $\Sigma_{g,0}$.
See, \emph{e.g.},~\cite{Birman88,FarbMarg11Book}.

\begin{figure}[tbhp]
  \centering
  \includegraphics[scale=0.7]{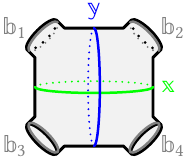}  
  \caption{Simple closed curves on the 4-punctured sphere
    $\Sigma_{0,4}$.}
  \label{fig:4-puncture}
\end{figure}

In the case of the 4-punctured sphere $\Sigma_{0,4}$,
the skein algebra is generated by~$\mathbb{x}$,~$\mathbb{y}$,
and~$\mathbb{b}_i$ in Fig.~\ref{fig:4-puncture}.
It is known that
the type-$C^\vee C_1$ DAHA is isomorphic to the skein
algebra on the 4-punctured sphere~\cite{Oblom04a,KHikami19a,BerestSamuel16a}.
See also~\cite{Zheda91a,Koorn07a,Koorn08a,Terwil13a} from a point of
view of the algebraic structure of the Askey--Wilson polynomials.
The Askey--Wilson operator $\ch \mathsf{Y}$ and $\ch \mathsf{X}$
respectively corresponds to the curves $\mathbb{y}$ and $\mathbb{x}$,
while the 4 parameters $(t_0,t_1,t_2,t_3)$ denote the  boundary
curves
$(\mathbb{b}_1, \mathbb{b}_3, \mathbb{b}_2, \mathbb{b}_4)$.
See~\cite{KHikami19a} for detail.
In addition,
the generators of 
the $SL(2;\mathbb{Z})$ action~\cite{Chered05Book} on $H_{q,\mathbf{t}}$
\begin{align}
  \label{sigma_R_sphere}
  \sigma_R=
  \begin{psmallmatrix}
    1 & 1 \\
    0 & 1
  \end{psmallmatrix}
  &:
    \begin{pmatrix}
      \mathsf{T}_0\\
      \mathsf{T}_1\\
      \mathsf{X}
    \end{pmatrix}
    \mapsto
    \begin{pmatrix}
      q^{-\frac{1}{2}} \mathsf{X}\mathsf{T}_0^{-1}\\
      \mathsf{T}_1\\
      \mathsf{X}
    \end{pmatrix}
    ,
  \\
  \label{sigma_L_sphere}
  \sigma_L=
  \begin{psmallmatrix}
    1 & 0 \\
    1 & 1
  \end{psmallmatrix}
  &:
    \begin{pmatrix}
      \mathsf{T}_0\\
      \mathsf{T}_1\\
      \mathsf{X}
    \end{pmatrix}
    \mapsto
    \begin{pmatrix}
      \mathsf{T}_0\\
      \mathsf{T}_1\\
      q^{\frac{1}{2}} \mathsf{T}_0\mathsf{X}^{-1} \mathsf{T}_1^{-1}
    \end{pmatrix}
    ,
\end{align}
can be interpreted as the half Dehn twists
along~$\mathbb{x}$ and~$\mathbb{y}$  respectively.

\section{DAHA on Double Torus and Skein Algebra}
\label{sec:Sigma20}

\subsection{\mathversion{bold}Skein Algebra
  and Mapping Class Group
  on \texorpdfstring{$\Sigma_{2,0}$}{Sigma 2,0}}

\begin{figure}[tbp]
  \centering
  \includegraphics[scale=0.7]{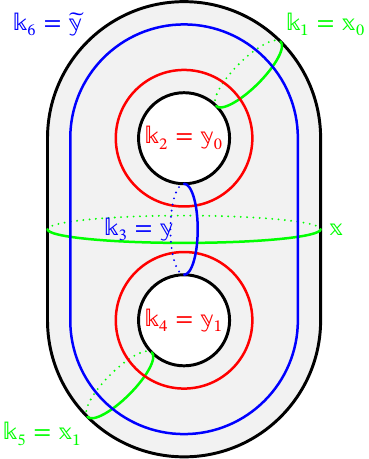}
  \caption{Simple closed curves on the double torus $\Sigma_{2,0}$.
}
  \label{fig:genus2surface}
\end{figure}

\begin{figure}[htbp]
  \centering
  \includegraphics[scale=0.7]{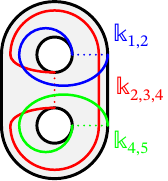}
  \qquad
  \includegraphics[scale=0.7]{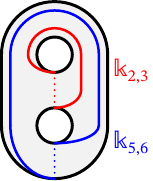}
  \qquad
  \includegraphics[scale=0.7]{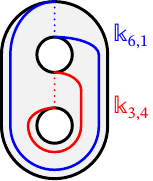}
  \qquad
  \includegraphics[scale=0.7]{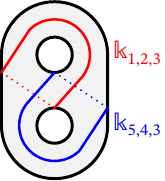}
  \qquad
  \includegraphics[scale=0.7]{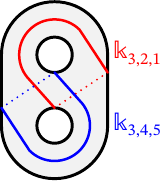}
  \caption{Several simple closed curves on $\Sigma_{2,0}$
    generated by the Dehn twists
    $\mathscr{D}_j$ from $\mathbb{k}_i$.}
  \label{fig:other_curves}
\end{figure}

The skein algebra $\Sk_A(\Sigma_{2,0})$ is generated by
$\mathbb{k}_1=\mathbb{x}_0$,
$\mathbb{k}_2=\mathbb{y}_0$,
$\mathbb{k}_3=\mathbb{y}$,
$\mathbb{k}_4=\mathbb{y}_1$,
and
$\mathbb{k}_5=\mathbb{x}_1$,
where
we label each simple closed curve on $\Sigma_{2,0}$ as in
Fig.~\ref{fig:genus2surface} following~\cite{KHikami19a}.
The Humphries generators  of the mapping class group are
$\mathscr{D}_i=\mathscr{D}_{\mathbb{k}_i}$
for $1\leq i \leq 5$,
and the mapping class group is
(see \emph{e.g.}~\cite{Birman88,Wajn83})
\begin{equation}
  \label{MCG_Sigma20}
  \begin{aligned}[b]
    \Mod(\Sigma_{2,0})
    & =
      \left\langle
      \mathscr{D}_1,\dots,
      \mathscr{D}_5
      ~\middle\vert~
      \begin{matrix}
        \text{$
        \mathscr{D}_{i,i+1,i}=\mathscr{D}_{i+1,i,i+1}
        $ for $1\leq i\leq 4$
        }
        \\
        \text{$
        \mathscr{D}_{i,j}=\mathscr{D}_{j,i}
        $ for $|i-j|>1$}
        \\
        \left(
          \mathscr{D}_{1,2,3,4,5}
        \right)^6=1
        \\
        \left(
        \mathscr{D}_{5,4,3,2,1,1,2,3,4,5}
        \right)^2=1
      \end{matrix}
      \right\rangle ,
  \end{aligned}
\end{equation}
where we mean
$\mathscr{D}_{i,\dots, j,k}=
\mathscr{D}_i
\dots \mathscr{D}_j\mathscr{D}_k$.
We note that the 2- and 3-chain relations~\cite{FarbMarg11Book} respectively give
\begin{gather}
  \label{2-chain}
  \left(\mathscr{D}_{1,2}\right)^6=
  \mathscr{D}_{\mathbb{x}} ,
  \\
  \label{3-chain}
  \left(\mathscr{D}_{1,2,3}\right)^4=\mathscr{D}_5^2   .
\end{gather}
See Fig.~\ref{fig:other_curves} for several simple closed curves
generated by the Dehn twists $\mathscr{D}_i$.
We mean for simplicity
$\mathbb{k}_{i,\pm j,\dots,\pm k}=
\left(\mathscr{D}_k^{\pm 1} \dots \mathscr{D}_j^{\pm 1}
\right)(\mathbb{k}_i)$.
These were  used in~\cite{Arthamo23a} in studies
of the character variety.

\subsection{$q$-Difference Operators}
In~\cite{KHikami19a} studied was the map
\begin{equation}
  \label{map_A}
  \mathcal{A}: \Sk_{A=q^{-\frac{1}{4}}}(\Sigma_{2,0})
  \to
  \End
  \mathbb{C}\left(q^{\frac{1}{4}},x_0,x_1\right)
  \left[ x+x^{-1} \right] .
\end{equation}
Therein
given are
for the curves in Fig.~\ref{fig:genus2surface} as
\begin{align}
  \label{A_for_xb}
  & \mathcal{A}(\mathbb{x}_b)
    = \ch(x_b),
  \\
  & \mathcal{A}(\mathbb{y}_b)
    =
    \I \,  q^{-\frac{1}{4}} G_0(x_b,x),
  \\
  & \mathcal{A}(\mathbb{x})
    = 
    \ch(x) ,
  \\
  \label{A_for_y}
  &
    \mathcal{A}(\mathbb{y})
    =
    \sum_{\epsilon=\pm} \omega(x^{\epsilon}) \left\{
    -   x^{-\epsilon} \left(x_0+\frac{q^{\frac{1}{2}}x^{\epsilon}}{x_0}\right)
    \left(x_1+\frac{q^{\frac{1}{2}}x^{\epsilon}}{x_1}\right) \, \eth^\epsilon
    +q^{\frac{1}{2}} \ch (x_0) \ch(x_1)
    \right\} ,
  \\
  %
  \label{A_tilde_y}
&
  \mathcal{A}(\widetilde{\mathbb{y}})
  =
  \sum_{\epsilon=\pm}
  \omega(x^{\epsilon}) \left\{
  K_0(x_0,x^{\epsilon})  \, K_0(x_1,x^{\epsilon}) \, \eth^{\epsilon}
  -G_0(x_0,x) \, G_0(x_1,x)
  \right\} ,
\end{align}
where  $b=0, 1$.
Here we have used
\begin{equation}
  \label{eq:22}
    \omega(x) =
  \frac{x \left(1+q^{\frac{1}{2}}x \right)}{
    q^{\frac{1}{2}}
    \left(1-x^2 \right)
    \left( 1- q^{\frac{1}{2}} x \right)
  } ,
\end{equation}
and
the $q$-difference operators $K_n(x_b,x)$ and $G_n(x_b,x)$ for
$n\in\mathbb{Z}$
are
defined by
\begin{align}
  \label{def_K_op}
  K_n(x_b,x)
  &=
    \frac{-x_b^{-n}}{1-x_b^2} \, \eth_b
    +
    \frac{
    x_b^n
    \left(q^{\frac{1}{2}} x + x_b^{2}\right)
    \left(q^{\frac{3}{2}} x+x_b^{2} \right)}{
    q\, x  \left( 1-x_b^{2} \right)} \, \eth_b^{-1} ,
  \\
  \label{def_G_op}
  G_n(x_b,x)
  &=
    \frac{-x_b^{-n}}{1-x_b^{2}} \, \eth_b
    +
    \frac{
    x_b^n
    \left(q^{\frac{1}{2}} x + x_b^{2}\right)
    \left(q^{\frac{1}{2}} +  x \, x_b^{2} \right)}{
        q^{\frac{1}{2}} \, x  \left( 1-x_b^{2} \right)}  \,
    \eth_b^{-1} ,
\end{align}
where the $q$-shift operators~$\eth_b$ for $b=0, 1$
are
\begin{equation}
  \label{eq:8}
  \left( \eth_0 f \right)(x,x_0,x_1)=f(x,  q^{\frac{1}{2}} x_0, x_1)    ,
  \qquad\qquad
  \left( \eth_1 f \right)(x,x_0,x_1)=f(x, x_0, q^{\frac{1}{2}} x_1) .
\end{equation}
Note that the $q$-difference operators $K_0(x_b, x)$,
$K_0(x_b,x^{-1})$, and
$G_0(x_b, x)$
are respectively
related to the raising operator, lowering operator, and the
eigen-operator for the type-$A_1$ Macdonald polynomials
a.k.a. the Rogers ultra-spherical polynomials~\cite{KHikami19a}.
We see that they  fulfill the following;
\begin{gather}
  \label{symmetric_G}
  G_n(x_b, x^{-1})=G_n(x_b, x) ,
  \\
  K_n(x_b,x) \, G_n(x_b, q\, x)
  = G_n(x_b, x) \, K_n(x_b, x) ,
  \\
  K_n(x_b, x^{-1}) \, K_n(x_b, q^{-1}x )
  - \left[ G_n(x_b, x) \right]^2
  = -q^{\frac{n}{2}}
  x^{-1}
  \left( q^{\frac{1}{2}} - x \right)^2 ,
  \\
  \label{relations_G_K_4}
  \begin{multlined}[b][.8\textwidth]
    q^{\frac{1}{2}}x \, (1-x^2) \left(
      G_n(x_b,x) \, G_n(x_b,q^{-1}x) - K_n(x_b,x) \, K_n(x_b,x^{-1})
    \right)
    \\
    +(1-q) \, x^2\left(
      K_n(x_b,x) - K_n(x_b,x^{-1})
    \right) G_n(x_b,q^{-1}x)
    \\
    =
    q^{\frac{n}{2}}\left(q^{\frac{1}{2}}-x \right)
    \left( 1- q^{\frac{1}{2}} x \right)
    \left( 1-x^2 \right) \frac{q-x_b^2}{1-x_b^2} ,
  \end{multlined}
  \\
  \label{G_K_recursion}
  \begin{pmatrix}
    K_{n+1}(x_b, x^{-1})
    \\
    G_{n+1}(x_b, q^{-1}x)
  \end{pmatrix}
  =
  \frac{1}{ x-q^{\frac{1}{2}}}
  \begin{pmatrix}
    x \ch(x_b) & - \frac{q^{\frac{1}{2}}+x\,x_b^2}{x_b}
    \\
    \frac{x+q^{\frac{1}{2}}x_b^2}{x_b} & -q^{\frac{1}{2}} \ch(x_b)
  \end{pmatrix}
  \begin{pmatrix}
    K_n(x_b,x^{-1})
    \\
    G_n(x_b,q^{-1}x)
  \end{pmatrix} .
\end{gather}

We note that
\begin{gather}
  \label{eq:4}
  \mathcal{A}(\mathbb{k}_{2,1^n})
  = \I \, q^{-\frac{n+1}{4}} G_n(x_0,x) ,
\end{gather}
which follows from
the skein algebra
\begin{equation}
  \label{eq:16}
  \mathbb{k}_1 \, \mathbb{k}_{2,1^n}=
  A \, \mathbb{k}_{2,1^{n-1}}
  +  A^{-1} \, \mathbb{k}_{2,1^{n+1}} ,
  \qquad
  \mathbb{k}_{2,1^n}
  =
  \begin{cases}
    \mathbb{k}_{2,\underbrace{1,\dots,1}_n},
    & \text{for $n\geq 0$},
    \\
    \mathbb{k}_{2,\underbrace{-1,\dots,-1}_{|n|}},
    & \text{for $n \leq 0$}.
  \end{cases}
\end{equation}


\subsection{Specialization of type-\texorpdfstring{$C^\vee C_1$}{CC1} DAHA}
Hereafter
we fix the parameters of the $C^\vee C_1$ DAHA as
$\mathbf{t}=\mathbf{t}_\star$,
\begin{equation}
  \label{special_t}
  \mathbf{t}_\star=
  \left(
    \I  \, x_0, \I \,  q^{-\frac{1}{2}} x_1,
    \I \,  x_0, \I \, x_1
  \right) ,
\end{equation}
to identify
the
curve~$\mathbb{b}_1$ (resp.~$\mathbb{b}_3$) with~$\mathbb{b}_2$
(resp.~$\mathbb{b}_4$) in Fig.~\ref{fig:4-puncture}.
At $\mathbf{t}_\star$~\eqref{special_t}
the Hecke operators~\eqref{AW_poly_rep} are read as
\begin{equation}
  \label{T_genus2}
  \begin{aligned}[]
    &\mathsf{T}_0
      \mapsto
      \I \, \frac{x}{q^{\frac{1}{2}} - x}
      \left(
      - \frac{q^{\frac{1}{2}}+ x\, x_0^2}{x\, x_0} \,
      \mathsf{s} \, \eth
      +
      x_0 + x_0^{~-1}
      \right) ,
    \\
    & 
      \mathsf{T}_1   \mapsto
      \I \left(
      \frac{1+q^{\frac{1}{2}}x}{q^{\frac{1}{2}}
      \left( 1-x^2 \right)}
      \frac{q^{\frac{1}{2}}x+x_1^2}{x_1} \,
      (\mathsf{s}-1)
      - q^{\frac{1}{2}} x_1^{-1}
      \right),
    \\
    & \mathsf{T}_0^\vee
      \mapsto
      q^{-\frac{1}{2}} {\mathsf{T}_0}^{-1} x
      ,
    \\
    & \mathsf{T}_1^\vee
     \mapsto
      x^{-1} {\mathsf{T}_1}^{-1}
      ,
  \end{aligned}
\end{equation}
which satisfy
the   Hecke relations
\begin{equation}
  \label{T01_Hecke}
  \begin{alignedat}{2}
    &\mathsf{T}_0-\mathsf{T}_0^{-1}
      =-\I \ch \left( x_0\right),
    &\qquad\qquad\qquad
    &\mathsf{T}_0^\vee-\left(\mathsf{T}_0^\vee\right)^{-1}
      =-\I \ch(x_0),
    \\
    &\mathsf{T}_1-\mathsf{T}_1^{-1}
      =-\I \ch\left(  q^{-\frac{1}{2}} \, x_1 \right),
    &
    &\mathsf{T}_1^\vee-\left(\mathsf{T}_1^\vee\right)^{-1}
      =-\I \ch(x_1) .
  \end{alignedat}
\end{equation}
The idempotent~\eqref{idempotent_AW} becomes
\begin{equation}
  \label{idempotent_e}
  \mathsf{e}
  \mapsto
  \left(1+\mathsf{s}\right)
  \frac{
    \left(q^{\frac{1}{2}}+x \right)
    \left( q^{\frac{1}{2}}+x \, x_1^2    \right)}{
    \left(1-x^2 \right) \left( q- x_1^2 \right)
  } .
\end{equation}
Note that
\begin{equation}
  \label{T1_e}
  \mathsf{T}_1 \, \mathsf{e}
  = -\I \, q^{\frac{1}{2}} \, x_1^{-1}
  \mathsf{e}
  = \mathsf{e} \, \mathsf{T}_1 ,
  \qquad\qquad
  \mathsf{T}_1^{-1} \, \mathsf{e} = \I \, q^{-\frac{1}{2}}\, x_1
  = \mathsf{e} \, \mathsf{T}_1^{-1}.
\end{equation}
and
\begin{equation}
  \label{x_T1_x}
  \mathsf{X} \mathsf{T}_1 \left( 1+ q^{\frac{1}{2}} \mathsf{X}\right)
  \,  \mathsf{e}
  =
  q^{\frac{1}{2}} \left( 1+q^{\frac{1}{2}}\mathsf{X}\right) \,
  \mathsf{T}_1^{-1} \, \mathsf{e} .
\end{equation}

The DAHA
at $\mathbf{t}_\star$~\eqref{special_t}
was employed so that the Askey--Wilson operator gives~\eqref{A_for_y}
as
\begin{equation}
  \label{A_y_TT}
  \mathcal{A}(\mathbb{y}) \,  \mathsf{e}
  = \ch(\mathsf{T}_1 \mathsf{T}_0) \, \mathsf{e} .
\end{equation}
Namely the Askey--Wilson polynomial~\eqref{AW_polynomial} is the eigen-polynomial of~\eqref{A_for_y},
\begin{equation*}
  \mathcal{A}(\mathbb{y}) \, P_m(x;q, \mathbf{t}_\star)
  = - \ch \left(\frac{q^{m+\frac{1}{2}}}{x_0 x_1} \right) \,
  P_m(x;q, \mathbf{t}_\star) .
\end{equation*}
It should be remarked that the  operator
$\mathcal{A}(\widetilde{\mathbb{y}})$, commuting with
the Askey--Wilson operator
as
$\mathcal{A}(\mathbb{y}) \, \mathcal{A}(\widetilde{\mathbb{y}})
=
\mathcal{A}(\widetilde{\mathbb{y}}) \, \mathcal{A}(\mathbb{y})$,
satisfies~\cite{KHikami24b}
\begin{multline*}
  \mathcal{A}(\widetilde{\mathbb{y}}) \,
  P_m(x;q, \mathbf{t}_\star)
  =
  \frac{\left(
      q^{m+2}- x_0^{2} x_1^{2} \right)^2}{
    q^{m+\frac{5}{2}} \left( 1-x_0^{ 2 }\right) \left( 1-x_1^{2} \right)}
  P_{m+1}(x ; q, \mathbf{t}_\star)
  \\
  -
  \frac{
    \left(q^{m+1}-x_0^{2}\right)^2 +\left(q^{m+1}-x_1^{2}\right)^2 
  }{
    q^{m+\frac{3}{2}} \left( 1-x_0^{2} \right) \left(1-x_1^{2} \right)}
  P_m(x; q, \mathbf{t}_\star)
  +
  \frac{(1-q^m)^2}{
    q^{m+\frac{1}{2}} \left( 1-x_0^{2} \right)\left(1-x_1^{2} \right)}
  P_{m-1}(x;q, \mathbf{t}_\star) .
\end{multline*}

\subsection{Heegaard Dual of Hecke Operators}
Our purpose is to rewrite
the map~\eqref{map_A}
in terms of the Iwahori--Hecke operators.
The motivation is  based on that
$S^3$ has a Heegaard splitting $S^3=H_1 \cup_{\Sigma_{2,0}} H_2$,
where $H_i$ is a 2-handlebody and~$\Sigma_{2,0}=\partial H_i$.
In gluing,
the meridians~$\mathbb{x}_b$ on~$H_1$
are  mapped to  the longitudes~$\mathbb{y}_b$
on~$H_2$,
and~$\mathbb{x}$ and~$\mathbb{y}$ are   to
$\mathbb{x}$ and $\widetilde{\mathbb{y}}$ respectively.
The fact that
$\mathbb{y}$ corresponds to the Askey--Wilson operator~\eqref{A_y_TT}
suggests
that there may exist
a ``Heegaard dual''
$\mathsf{U}_0$ and $\mathsf{U}_1$
of the Hecke operators $\mathsf{T}_0$ and
$\mathsf{T}_1$~\eqref{T_genus2} for~$\widetilde{\mathbb{y}}$.

\begin{definition}
  \label{def:U_operator}
We define the  representation of $\mathsf{U}_0$ and
$\mathsf{U}_1$ by
\begin{align}
  \label{def_00}
  & \mathsf{U}_0
   \mapsto
    \frac{q^{-\frac{1}{4}}x}{q^{\frac{1}{2}}-x} \,
    K_0(x_0, x^{-1}) \,
    \mathsf{s} \, \eth
    -
    \frac{q^{-\frac{1}{4}}x}{q^{\frac{1}{2}}-x} \, G_0(x_0,x) ,
  \\
  \label{def_01}
  &
    \mathsf{U}_1
   \mapsto
        -\frac{x \left( 1+q^{\frac{1}{2}}x \right)}{
    q^{\frac{1}{4}} \left( 1-x^2 \right)} \,
    K_0(x_1, x) \,
        (\mathsf{s}-1)
        +\frac{q^{\frac{1}{4}}}{1-q^{\frac{1}{2}}x}
        \left( G_0(x_1,x) - q^{\frac{1}{2}} x \, K_0(x_1,x) \right)  ,
  \\
  &  \mathsf{U}_0^\vee
   \mapsto
    q^{-\frac{1}{2}} \mathsf{U}_0^{-1} x ,
  \\
  & \mathsf{U}_1^\vee
   \mapsto
    x^{-1} \mathsf{U}_1^{-1} ,
\end{align}
where
$K_n(x_b,x)$ and  $G_n(x_b,x)$ are given  in~\eqref{def_K_op}
and~\eqref{def_G_op}.
\end{definition}
The invertibilities of $\mathsf{U}_0$
and $\mathsf{U}_1$
can be
checked 
using~\eqref{symmetric_G}--\eqref{relations_G_K_4}
by
\begin{align}
  \label{eq:11}
  &\mathsf{U}_0^{-1}
   \mapsto
    \frac{q^{-\frac{1}{4}}x}{q^{\frac{1}{2}}-x}  \, K_0(x_0, x^{-1}) \,
    \mathsf{s} \, \eth
    -
    \frac{q^{\frac{1}{4}}}{q^{\frac{1}{2}}-x} \, G_0(x_0,x)
    ,
  \\
  &
    \mathsf{U}_1^{-1}
    \mapsto
        \left\{
        (\mathsf{s}+1)
        \frac{q^{\frac{1}{2}}+x}{q^{\frac{1}{4}}(1-x^2)} \,
        K_0(x_1, x^{-1})
        +\frac{q^{\frac{1}{4}}}{q^{\frac{1}{2}}-x}
        \left(
        x\, G_0(x_1, q^{-1}x)
        - q^{\frac{1}{2}} \, K_0(x_1,x^{-1})
        \right)
        \right\}
        \frac{1-x_1^2}{q-x_1^2}  .
\end{align}
By construction, we have
an analogue of~\eqref{prod_TTTT}
\begin{equation}
  \label{prod_UUUU}
  \mathsf{U}_1^\vee  \mathsf{U}_1  \mathsf{U}_0  \mathsf{U}_0^\vee
  =q^{-\frac{1}{2}}.
\end{equation}
Moreover we get the followings.
\begin{lemma}
  \label{lemma:dual_condition}
  \begin{enumerate}
  \item 
    \begin{gather}
      \label{prod_UTUT}
      \mathsf{U}_0  \mathsf{T}_0 \mathsf{U}_0^\vee
      = - q^{\frac{1}{2}}  \,\mathsf{T}_0 ,
    \end{gather}

  \item
    \begin{equation}
      \label{relation_UTXUT}
      \mathsf{U}_1 \mathsf{T}_1^{-1}  \mathsf{U}_1^\vee
      \mathsf{T}_1
      \mapsto
      -q
      - \frac{(1-q)
        \left( 1+q^{\frac{1}{2}}x \right)
        \left(q^{\frac{1}{2}}x+x_1^2\right)
      }{\left(1-x^2 \right) \left(q-x_1^2\right)}
      \left( \mathsf{s} - 1 \right) .
    \end{equation}
  \end{enumerate}
\end{lemma}

The representations in Def.~\ref{def:U_operator}
and the inverses
show that
the Heegaard dual operators
satisfy the
Hecke-type relations.
We have for $\mathsf{U}_0$ and $\mathsf{U}_0^\vee$
\begin{align}
  \label{U0_Hecke}
  &
    \begin{aligned}
      \mathsf{U}_0 - \mathsf{U}_0^{~ -1}
      & =
        q^{-\frac{1}{4}} \, G_0(x_0,x),
      \\
      \mathsf{U}^\vee_0 - \left( \mathsf{U}^\vee_0 \right)^{-1}
      & =
        q^{-\frac{1}{4}} \, G_0(x_0,x) .
    \end{aligned}
\end{align}
For $\mathsf{U}_1$ and
$\mathsf{U}_1^\vee$,
we
have preferable  expressions
with the symmetrizer $\mathsf{e}$~\eqref{idempotent_e}
\begin{equation}
  \label{U1_Hecke_e}
  \begin{aligned}
    \left(
    q^{-\frac{1}{2}} \mathsf{U}_1 - q^{\frac{1}{2}}\mathsf{U}_1^{-1}
    \right) \,\mathsf{e}
    & =
      q^{-\frac{1}{4}} G_0(x_1,x) \, \mathsf{e} ,
    \\
    \left(
    \mathsf{U}_1^\vee - \left(\mathsf{U}_1^\vee \right)^{-1}
    \right) \, \mathsf{e}
    & =
      q^{-\frac{1}{4}} G_0(x_1,x) \, \mathsf{e} ,
  \end{aligned}
\end{equation}
which  can be seen by use of
an analogous identity to~\eqref{x_T1_x},
\begin{equation}
  \label{eq:20}
  \mathsf{X} \mathsf{U}_1 \left( 1+ q^{\frac{1}{2}} \mathsf{X}\right)
  \,  \mathsf{e}
  =
  q^{\frac{1}{2}} \left( 1+q^{\frac{1}{2}}\mathsf{X}\right) \,
  \mathsf{U}_1^{-1} \, \mathsf{e} .
\end{equation}

Furthermore we can prove
\begin{equation}
  \label{eq:9}
  \mathcal{A}(\widetilde{\mathbb{y}}) \,\mathsf{e}
  =
  \ch\left(
    \mathsf{U}_1  \mathsf{U}_0
  \right)  \mathsf{e}.
\end{equation}
In summary,
all the generators $\mathbb{k}_i$ for
the skein algebra
$\Sk_{A=q^{-\frac{1}{4}}}(\Sigma_{2,0})$ given in~\eqref{A_for_xb}--\eqref{A_tilde_y}
can be written as follows.
\begin{prop}
\begin{equation}
  \label{A_from_DAHA}
  \begin{aligned}
    & \mathcal{A}(\mathbb{k}_1)
      =
      \mathcal{A}(\mathbb{x}_0)  =
      \ch \left( \I \, \mathsf{T}_0\right)
      =
      \ch \left( \I \, \mathsf{T}_0^\vee \right) ,
    \\
    & \mathcal{A}(\mathbb{k}_2)
      =
      \mathcal{A}(\mathbb{y}_0)
      =
      \ch \left( \I \, \mathsf{U}_0 \right)
      =
      \ch \left( \I \, \mathsf{U}_0^\vee \right) ,
      \\
    & \mathcal{A}(\mathbb{k}_3) \, \mathsf{e}
      =
      \mathcal{A}(\mathbb{y})  \, \mathsf{e}
      =
      \ch (\mathsf{T}_1 \mathsf{T}_0)  \, \mathsf{e}
      =
      \ch (\mathsf{T}_0 \mathsf{T}_1)  \, \mathsf{e} ,
    \\
    & \mathcal{A}(\mathbb{k}_4) \, \mathsf{e}
      =
      \mathcal{A}(\mathbb{y}_1)\,\mathsf{e}
      =
      \ch \left( \I \, q^{-\frac{1}{2}} \mathsf{U}_1 \right)  \mathsf{e}
      =
      \ch\left( \I \, \mathsf{U}_1^\vee \right)  \mathsf{e},
    \\
    & \mathcal{A}(\mathbb{k}_5) \, \mathsf{e}
    =
      \mathcal{A}(\mathbb{x}_1) \, \mathsf{e}
    =
    \ch \left( \I \, q^{-\frac{1}{2}} \mathsf{T}_1 \right)  \mathsf{e}
    =
    \ch \left( \I \, \mathsf{T}_1^\vee \right)  \mathsf{e},
    \\
    & \mathcal{A}(\mathbb{k}_6) \, \mathsf{e}
      =
      \mathcal{A}(\widetilde{\mathbb{y}}) \,\mathsf{e}
      =
      \ch\left(
      \mathsf{U}_1  \mathsf{U}_0
      \right)  \mathsf{e} .
  \end{aligned}
\end{equation}
\end{prop}
We can thus regard the map~\eqref{map_A}
as
\begin{equation}
  \label{spherical_A}
  \mathcal{A}: \Sk_{A=q^{-\frac{1}{4}}}(\Sigma_{2,0})
  \to
  SH_{q,\mathbf{t}_\star}^{gen}
\end{equation}
where
$SH_{q,\mathbf{t}_\star}^{gen}$
is
a spherical subalgebra of
our generalized DAHA,
\begin{equation}
  \label{general_H_q}
  H_{q,\mathbf{t}_\star}^{gen}
  =
  \left\langle \mathsf{T}_0^{\pm 1},
    \mathsf{T}_1^{\pm 1}, \mathsf{X}^{\pm 1},
    \mathsf{U}_0^{\pm 1}, \mathsf{U}_1^{\pm 1}
    ~\middle\vert~
    \begin{matrix}
      \text{the Hecke
      relations~\eqref{T01_Hecke},~\eqref{U0_Hecke},~\eqref{U1_Hecke_e}}
      \\
      \mathsf{T}_0 \mathsf{U}_0^{-1} \mathsf{X} \mathsf{T}_0^{-1}
      \mathsf{U}_0
      =-q
      \\
      \mathsf{U}_1 \mathsf{T}_1^{-1} \mathsf{X}^{-1}
      \mathsf{U}_1^{-1}\mathsf{T}_1 \,
      \mathsf{e}
      = -q \,
      \mathsf{e}
    \end{matrix}
  \right\rangle .
\end{equation}
The conditions are from Lemma~\ref{lemma:dual_condition}.

\subsection{Automorphisms}

By definition of
$H_{q,\mathbf{t}_\star}^{gen}$~\eqref{general_H_q},
we find  the automorphisms
$\mathscr{T}_i=\mathscr{T}_{\mathbb{k}_i}$ as follows.

\begin{prop}
  We have the automorphisms of $H_{q,\mathbf{t}_\star}^{gen}$~\eqref{general_H_q};
\begin{align}
  \label{auto_T1}
  \mathscr{T}_1=
  \mathscr{T}_{\mathbb{x}_0}:
  &
    \begin{pmatrix}
      \mathsf{T}_0
      \\
      \mathsf{T}_1
      \\
      \mathsf{X}
      \\
      \mathsf{U}_0
      \\
      \mathsf{U}_1
    \end{pmatrix}
    \mapsto
    \begin{pmatrix}
      \mathsf{T}_0
      \\
      \mathsf{T}_1
      \\
      \mathsf{X}
      \\
      -\I \, q^{\frac{1}{4}}\mathsf{U}_0 \mathsf{T}_0^{-1}
      \\
      \mathsf{U}_1
    \end{pmatrix}
    ,
  \\
  \mathscr{T}_2=
  \mathscr{T}_{\mathbb{y}_0}:
  &
    \begin{pmatrix}
      \mathsf{T}_0
      \\
      \mathsf{T}_1
      \\
      \mathsf{X}
      \\
      \mathsf{U}_0
      \\
      \mathsf{U}_1
    \end{pmatrix}
    \mapsto
    \begin{pmatrix}
      \I \, q^{-\frac{1}{4}}\mathsf{U}_0 \mathsf{T}_0
      \\
      \mathsf{T}_1
      \\
      \mathsf{X}
      \\
      \mathsf{U}_0
      \\
      \mathsf{U}_1
    \end{pmatrix}
    ,
  \\
  \label{auto_T3}
  \mathscr{T}_3=
  \mathscr{T}_{\mathbb{y}}:
  &
    \begin{pmatrix}
      \mathsf{T}_0
      \\
      \mathsf{T}_1
      \\
      \mathsf{X}
      \\
      \mathsf{U}_0
      \\
      \mathsf{U}_1
    \end{pmatrix}
    \mapsto
    \begin{pmatrix}
      \mathsf{T}_0
      \\
      \mathsf{T}_1
      \\
      \left( \mathsf{T}_0 \mathsf{T}_1 \right)^{-1}
      \mathsf{X}
      \mathsf{T}_1 \mathsf{T}_0
      \\
      q^{-\frac{1}{4}} \left( \mathsf{T}_0 \mathsf{T}_1\right)^{-1}
      \mathsf{U}_0
      \\
      q^{\frac{1}{4}}
      \mathsf{U}_1
      \mathsf{T}_0 \mathsf{T}_1
    \end{pmatrix}
    ,
  \\
  \mathscr{T}_4=
  \mathscr{T}_{\mathbb{y}_1}:
  &
    \begin{pmatrix}
      \mathsf{T}_0
      \\
      \mathsf{T}_1
      \\
      \mathsf{X}
      \\
      \mathsf{U}_0
      \\
      \mathsf{U}_1
    \end{pmatrix}
    \mapsto
    \begin{pmatrix}
      \mathsf{T}_0
      \\
      \I \, q^{-\frac{1}{4}}
      \left( \mathsf{U}_1\mathsf{X}\right)^{-1}
      \mathsf{T}_1
      \\
      \mathsf{X}
      \\
      \mathsf{U}_0
      \\
      \mathsf{U}_1
    \end{pmatrix}
    ,
  \\
  \label{auto_T5}
  \mathscr{T}_5=
  \mathscr{T}_{\mathbb{x}_1}:
  &
    \begin{pmatrix}
      \mathsf{T}_0
      \\
      \mathsf{T}_1
      \\
      \mathsf{X}
      \\
      \mathsf{U}_0
      \\
      \mathsf{U}_1
    \end{pmatrix}
    \mapsto
    \begin{pmatrix}
      \mathsf{T}_0
      \\
      \mathsf{T}_1
      \\
      \mathsf{X}
      \\
      \mathsf{U}_0
      \\
      -\I \, q^{\frac{1}{4}}
      \mathsf{U}_1
      \mathsf{X} \mathsf{T}_1
    \end{pmatrix}
    .
\end{align}
\end{prop}

We note that the map $\mathscr{T}_3$~\eqref{auto_T3} originates from
the  Dehn
twist~$\sigma_L^{-2}$ on $\Sigma_{0,4}$~\eqref{sigma_L_sphere}.

Our claim is as follows.
\begin{prop}
  We have a commutative diagram;
  \begin{equation*}
    \begin{tikzcd}
      \Sk_{A=q^{-\frac{1}{4}}}(\Sigma_{2,0})\arrow[r, "\mathscr{D}_i"]
      \arrow[d, "\mathcal{A}"']
      & \Sk_{A=q^{-\frac{1}{4}}}(\Sigma_{2,0})
      \arrow[d, "\mathcal{A}"]
      \\
      SH_{q,\mathbf{t}_\star}^{gen}
      \arrow[r, "\mathscr{T}_i"']
      &
      SH_{q,\mathbf{t}_\star}^{gen}
    \end{tikzcd}
  \end{equation*}
\end{prop}

We shall check the relations in~\eqref{MCG_Sigma20}.
From the  definitions~\eqref{auto_T1}--\eqref{auto_T5}
it is straightforward to see case-by-case that both the
braid relations and the commutativities hold;
\begin{gather}
  \mathscr{T}_{i,i+1,i}
  =
  \mathscr{T}_{i+1,i,i+1} ,
  \qquad
  \text{for $1\leq i \leq 4$},
  \\
  \mathscr{T}_{i,j}=\mathscr{T}_{j,i},
  \qquad
  \text{for $|i-j|>1$,}
\end{gather}
where we mean $\mathscr{T}_{i,\dots,j,k}=
\mathscr{T}_i
\dots
\mathscr{T}_j \mathscr{T}_k$.
We can also find that
\begin{align*}
    \mathscr{T}_{1,2,3,4,5}:
  &
    \begin{pmatrix}
      \mathsf{T}_0
      \\
      \mathsf{T}_1
      \\
      \mathsf{X}
      \\
      \mathsf{U}_0
      \\
      \mathsf{U}_1
    \end{pmatrix}
    \mapsto
    \begin{pmatrix}
      \mathsf{U}_0
      \\
      \I \, q^{-\frac{1}{2}} \mathsf{U}_0^{-1} \mathsf{T}_1^{-1}
      \mathsf{X}^{-1} \mathsf{U}_1^{-1} \mathsf{T}_1
      \\
      \mathsf{T}_1^{-1} \mathsf{U}_0^{-1}
      \mathsf{X}
      \mathsf{T}_1\mathsf{U}_0
      \\
      -\I \,\mathsf{T}_1^{-1}
      \mathsf{T}_0^{-1}
      \\
      \mathsf{T}_1
    \end{pmatrix}
    ,
  \\
  \mathscr{T}_{5,4,3,2,1}:
  &
  \begin{pmatrix}
      \mathsf{T}_0
      \\
      \mathsf{T}_1
      \\
      \mathsf{X}
      \\
      \mathsf{U}_0
      \\
      \mathsf{U}_1
    \end{pmatrix}
    \mapsto
    \begin{pmatrix}
      \I \, q \mathsf{T}_1^{-1} \mathsf{U}_1 \mathsf{X}
      \mathsf{T}_1\mathsf{U}_0
      \\
      -q^{-\frac{1}{2}} \mathsf{X}^{-1}
      \left(\mathsf{T}_1^{-1}\mathsf{U}_1\mathsf{X} \mathsf{T}_1 \right)^{-1}
      \\
      \left( \mathsf{T}_1^{-1} \mathsf{U}_1 \mathsf{X} \mathsf{T}_1 \right)
      \mathsf{X} \mathsf{T}_0^{-1}
      \left( \mathsf{T}_1^{-1} \mathsf{U}_1 \mathsf{X} \mathsf{T}_1
      \right)^{-1}
      \mathsf{T}_0
      \\
      q^{-\frac{1}{2}}
      \left( \mathsf{T}_1^{-1} \mathsf{U}_1 \mathsf{X} \mathsf{T}_1
      \right)
      \mathsf{X}\mathsf{T}_0^{-1}
      \left( \mathsf{T}_1^{-1} \mathsf{U}_1 \mathsf{X} \mathsf{T}_1 \right)^{-1}
      \\
      \I \,
      \mathsf{U}_1 \mathsf{X} \mathsf{T}_1 \mathsf{T}_0
      \mathsf{X}^{-1}
      \left( \mathsf{T}_1^{-1} \mathsf{U}_1 \mathsf{X} \mathsf{T}_1 \right)^{-1}
    \end{pmatrix}
    ,
  \\
\end{align*}
which result in
\begin{equation}
  \label{T_12345-6}
  \left(\mathscr{T}_{1,2,3,4,5}\right)^6
  =
  \left(\mathscr{T}_{5,4,3,2,1,1,2,3,4,5}\right)^2:
  \begin{pmatrix}
    \mathsf{T}_0
    \\
    \mathsf{T}_1
    \\
    \mathsf{X}
    \\
    \mathsf{U}_0
    \\
    \mathsf{U}_1
  \end{pmatrix}
  \mapsto
  \left(
    \mathsf{T}_1^{-1} \mathsf{U}_1 \mathsf{X} \mathsf{T}_1
    \mathsf{U}_1^{-1}
  \right)
  \begin{pmatrix}
    \mathsf{T}_0
    \\
    \mathsf{T}_1
    \\
    \mathsf{X}
    \\
    \mathsf{U}_0
    \\
    \mathsf{U}_1
  \end{pmatrix}
  \left(
    \mathsf{T}_1^{-1} \mathsf{U}_1 \mathsf{X} \mathsf{T}_1
    \mathsf{U}_1^{-1}
  \right)^{-1} .
\end{equation}
As seen from~\eqref{relation_UTXUT},
the operator
$\mathsf{U}_1 \mathsf{T}_1^{-1}  \mathsf{X}^{-1}
\mathsf{U}_1^{-1} \mathsf{T}_1
$
acts as a scalar  on
the symmetric Laurent polynomials,
and the maps~\eqref{T_12345-6} are identities on~$\mathbb{C}(q^{\frac{1}{4}},x_0,x_1)[x+x^{-1}]$.

It should be noted that,
for the 3-chain relation~\eqref{3-chain}, we have
\begin{equation*}
  \left(\mathscr{T}_{1,2,3}\right)^4
  :
  \begin{pmatrix}
    \mathsf{T}_0 \\
    \mathsf{T}_1\\
    \mathsf{X}    \\
    \mathsf{U}_0 \\
    \mathsf{U}_1
  \end{pmatrix}
  \mapsto
  \begin{pmatrix}
    \mathsf{T}_1^{-1} \mathsf{T}_0 \mathsf{T}_1 \\
    \mathsf{T}_1\\
    \mathsf{T}_1^{-1} \mathsf{X} \mathsf{T}_1 \\
    \mathsf{T}_1^{-1} \mathsf{U}_0 \mathsf{T}_1 \\
    \mathsf{U}_1 \mathsf{X} \mathsf{T}_1^2
  \end{pmatrix}
  ,
  \qquad \qquad
  \left(\mathscr{T}_5\right)^2
  :
  \begin{pmatrix}
    \mathsf{T}_0 \\
    \mathsf{T}_1\\
    \mathsf{X}    \\
    \mathsf{U}_0 \\
    \mathsf{U}_1
  \end{pmatrix}
  \mapsto
  \begin{pmatrix}
    \mathsf{T}_0 \\
    \mathsf{T}_1\\
    \mathsf{X}    \\
    \mathsf{U}_0 \\
    -q^{\frac{1}{2}}\mathsf{U}_1 \mathsf{X} \mathsf{T}_1 \mathsf{X} \mathsf{T}_1
  \end{pmatrix}
  .
\end{equation*}
Both actions on the generators~\eqref{A_from_DAHA}  are same.

We explicitly
give the map for the curves
in Fig.~\ref{fig:other_curves}.
We have checked  the consistency with the skein algebra,
\emph{e.g.}~\eqref{example_skein},
of the curves.
\begin{align*}
  &
    \begin{aligned}
      \mathcal{A}(\mathbb{k}_{1,2})
      &=
        \mathscr{T}_2(\mathcal{A}(\mathbb{k}_1))
        =
        \ch\left( - q^{-\frac{1}{4}} \mathsf{U}_0 \mathsf{T}_0\right)
        =
        \ch\left(
        q^{\frac{1}{4}} \mathsf{X}^{-1} \mathsf{U}_0 \mathsf{T}_0
        \right)
      \\
      & =
        \I \, G_{-1}(x_0,x),
        %
    \end{aligned}
  \\
  &
    \begin{aligned}
      \mathcal{A}(\mathbb{k}_{2,3}) \, \mathsf{e} 
      & =
        \mathscr{T}_3(\mathcal{A}(\mathbb{k}_2)) \, \mathsf{e}
        =
        \ch\left( \I\, q^{-\frac{1}{4}} (\mathsf{T}_0 \mathsf{T}_1 )^{-1}
        \mathsf{U}_0 \right)  \mathsf{e}
        =
        \ch\left(
        -\I \, q^{\frac{1}{4}} (\mathsf{T}_1 \mathsf{T}_0)^{-1}
        \mathsf{X}^{-1} \mathsf{U}_0
        \right)  \mathsf{e}
        \\
      & =
        \I
        \sum_{\epsilon=\pm} \omega(x^{\epsilon})
        \left(
        -
        K_1(x_0,x^{\epsilon}) \,
        \frac{x_1^2+q^{\frac{1}{2}} x^{\epsilon}}{x_1}
        \,\eth^{\epsilon} + G_1(x_0,x) \,\ch(x_1) 
        \right)  \mathsf{e}
        ,
    \end{aligned}
  \\
  &
    \begin{aligned}
      \mathcal{A}(\mathbb{k}_{3,4}) \, \mathsf{e} 
      & =
        \mathscr{T}_3^{-1}(\mathcal{A}(\mathbb{k}_4)) \, \mathsf{e}
        =\ch \left( \I \, q^{-\frac{3}{4}}
        \mathsf{U}_1 ( \mathsf{T}_0 \mathsf{T}_1)^{-1}
        \right)  \mathsf{e}
        =
        \ch \left(
        -\I \, q^{-\frac{1}{4}}
        \mathsf{U}_1 \mathsf{X}
        ( \mathsf{T}_1 \mathsf{T}_0)^{-1}
        \right)  \mathsf{e}
      \\
      & =
        \I \sum_{\epsilon=\pm}
        \omega(x^{\epsilon})
        \left(
        -
        \frac{x^{-\epsilon} x_0^2+q^{\frac{1}{2}}}{x_0}
         K_{-1}(x_1,x^{\epsilon}) \, \eth^{\epsilon}
        +q^{\frac{1}{2}} \ch(x_0) \, G_{-1}(x_1,x)
        \right)  \mathsf{e}
        ,
    \end{aligned}
  \\
  &
    \begin{aligned}
      \mathcal{A}(\mathbb{k}_{4,5}) \,  \mathsf{e} 
      &=
        \mathscr{T}_5(\mathcal{A}(\mathbb{k}_4)) \,
        \mathsf{e}
        =
        \ch \left( q^{-\frac{1}{4}} \mathsf{U}_1\mathsf{X} \mathsf{T}_1 \right) 
        \mathsf{e}
        =
        \ch \left(
        -q^{\frac{1}{4}} \mathsf{U}_1 \mathsf{X} \mathsf{T}_1 \mathsf{X}
        \right)  \mathsf{e}
      \\
      & =
        \I \, q^{-\frac{1}{2}} G_1(x_1,x) \, \mathsf{e},
    \end{aligned}
  \\
  &
    \begin{aligned}
      \mathcal{A}(\mathbb{k}_{5,6}) \, \mathsf{e} 
      & =
        \mathscr{T}_5^{-1} \left(
        \mathcal{A}(\mathbb{k}_6)
        \right)  \mathsf{e}
        =
        \ch\left( \I \, q^{-\frac{1}{4}} \mathsf{U}_1
        (\mathsf{X} \mathsf{T}_1)^{-1} \mathsf{U}_0
        \right) \mathsf{e}
      \\
      & =
        q^{\frac{1}{4}}
        \sum_{\epsilon=\pm} \omega(x^{\epsilon})
        \left(  K_{0}(x_0,x^{\epsilon}) \, K_{-1}(x_1,x^{\epsilon})
        \, \eth^{\epsilon}
        - G_0(x_0,x) \, G_{-1}(x_1,x)
        \right) \mathsf{e}
        ,
    \end{aligned}
  \\
  &
    \begin{aligned}
      \mathcal{A}(\mathbb{k}_{6,1}) \, \mathsf{e} 
      &=
        \mathscr{T}_1 \left(
        \mathcal{A}(\mathbb{k}_6)
        \right)  \mathsf{e}
        =
        \ch\left(-\I \,  q^{\frac{1}{4}} \mathsf{U}_1 \mathsf{U}_0
        \mathsf{T}_0^{-1}
        \right)  \mathsf{e}
      \\
      & =
        q^{-\frac{1}{4}}
        \sum_{\epsilon=\pm}
        \omega(x^{\epsilon})
        \left(
        K_1(x_0,x^{\epsilon}) \, K_0(x_1,x^{\epsilon}) \, \eth^{\epsilon}
        - G_1(x_0,x)\, G_0(x_1,x)
        \right) \mathsf{e}
        ,
    \end{aligned}
  \\
  &
    \begin{aligned}
      \mathcal{A}(\mathbb{k}_{1,2,3}) \, \mathsf{e} 
      &=
        \mathscr{T}_3\left(
        \mathcal{A}(\mathbb{k}_{1,2})
        \right)  \mathsf{e}
        =
        \ch \left(
        q^{\frac{1}{2}} \mathsf{T}_1^{-1}\mathsf{X}^{-1}
        \mathsf{U}_0
        \right)  \mathsf{e}
        =
        \ch \left(
        (\mathsf{X}\mathsf{T}_1 \mathsf{T}_0)^{-1}
        \mathsf{U}_0 \mathsf{T}_0\right)  \mathsf{e}
      \\
      & =
        \I \, q^{\frac{1}{4}}
        \sum_{\epsilon=\pm}
        \omega(x^{\epsilon})
        \left(
        - K_0(x_0,x^{\epsilon}) \,
        \frac{x_1^2+q^{\frac{1}{2}}x^{\epsilon}}{x_1}
        \eth^{\epsilon}
        +G_0(x_0,x) \,\ch(x_1)
        \right)\mathsf{e}
        ,
    \end{aligned}
  \\
    &
    \begin{aligned}
      \mathcal{A}(\mathbb{k}_{5,4,3}) \, \mathsf{e}
      & =
        \left( \mathscr{T}_5^{-1} \mathscr{T}_3 \right)
        (\mathcal{A}(\mathbb{k}_4)) \, \mathsf{e}
        =
        \ch (\mathsf{U}_1 \mathsf{T}_0) \, \mathsf{e}
        =
        \ch\left(
        -q^{\frac{1}{2}} \mathsf{U}_1 \mathsf{T}_1^{-1}\mathsf{X}^{-1}
        \mathsf{T}_0\mathsf{T}_1
        \right) \, \mathsf{e}
      \\
      & =
        \I \, q^{\frac{1}{4}}
        \sum_{\epsilon=\pm}
        \omega(x^\epsilon)
        \left(
        -
        \frac{x_0^2+q^{\frac{1}{2}} x^\epsilon}{x_0}
        \, K_0(x_1,x^\epsilon) \, \eth^\epsilon
        +\ch(x_0) \, G_0(x_1,x)
        \right) \mathsf{e} ,
    \end{aligned}
  \\
  &
    \begin{aligned}
    \mathcal{A}(\mathbb{k}_{2,3,4}) \, \mathsf{e} 
      & =
        \mathscr{T}_2^{-1}\left(
        \mathcal{A}(\mathbb{k}_{3,4})
        \right) \, \mathsf{e}
        =
        \ch\left(
        -q^{-1} \mathsf{U}_1
        (\mathsf{T}_0 \mathsf{T}_1)^{-1} \mathsf{U}_0
        \right) \, \mathsf{e}
        =
        \ch\left(
        -q^{\frac{1}{2}} \mathsf{U}_1 \mathsf{U}_0
        (\mathsf{T}_1 \mathsf{T}_0)^{-1}
        \right) \, \mathsf{e} 
      \\
      & =
        \sum_{\epsilon=\pm}
        \omega(x^{\epsilon})
        \left(
        K_1(x_0,x^{\epsilon}) \, K_{-1}(x_1,x^{\epsilon}) \, \eth^{\epsilon}
        -
        G_1(x_0,x) \, G_{-1}(x_1,x)
        \right) \mathsf{e} 
        ,
    \end{aligned}
  \\
  &
    \begin{aligned}
    \mathcal{A}(\mathbb{k}_{3,4,5}) \, \mathsf{e} 
      &=
        \mathscr{T}_5 \left(
        \mathcal{A}(\mathbb{k}_{3,4})
        \right)  \mathsf{e}
        =
        \ch
        \left(
        q^{-\frac{1}{2}} \mathsf{U}_1\mathsf{X} \mathsf{T}_0^{-1}
        \right)  \mathsf{e}
        =
        \ch\left(
        - \mathsf{U}_1 \mathsf{X} \mathsf{T}_1 \mathsf{X}
        (\mathsf{T}_1 \mathsf{T}_0)^{-1}
        \right)  \mathsf{e}
      \\
      & =
        \I \, q^{\frac{1}{4}}
        \sum_{\epsilon=\pm}
        \omega(x^{\epsilon})
        \left(
        -
        \frac{q^{-\frac{1}{2}} x^{-\epsilon} x_0^2+1}{x_0} \,
        K_0(x_1,x^{\epsilon}) \,
        \eth^{\epsilon}
        +
        \ch(x_0) \,
        G_0(x_1,x)
        \right) \mathsf{e}
        ,
    \end{aligned}
  \\
    &
    \begin{aligned}
      \mathcal{A}(\mathbb{k}_{3,2,1})\, \mathsf{e}
      & =
        \left(\mathscr{T}_1 \mathscr{T}_2 \right)
        \left(
        \mathcal{A}(\mathbb{k}_3)
        \right)  \mathsf{e}
        =
        \ch(\mathsf{T}_1 \mathsf{U}_0) \, \mathsf{e}
        =
        \ch(\mathsf{U}_0 \mathsf{T}_1) \, \mathsf{e}
      \\
      & =
        \I \, q^{\frac{1}{4}}
        \sum_{\epsilon=\pm}
        \omega(x^\epsilon)
        \left(
        -
        K_0(x_0,x^\epsilon) \,
        \frac{
        q^{-\frac{1}{2}} x^{-\epsilon} x_1^2+1}{x_1}
        \, \eth^\epsilon
        +G_0(x_0,x) \, \ch(x_1)
        \right) \mathsf{e} .
    \end{aligned}
\end{align*}
Here we avoid to use $\mathscr{T}_4$ due to that $\mathsf{T}_1$,
used
for the idempotent~$\mathsf{e}$, is no
longer  invariant.

\subsection{Rational Tangles}
In \cite{JConway70a},
Conway introduced tangle operations, and showed that
continued fraction can be assigned to a certain family of knots and links.
In view from  links on the  double torus~$\Sigma_{2,0}$, the tangle
operations correspond
to the Dehn
twists along~$\mathbb{x}$ and~$\mathbb{y}$ acting
on the  curve~$\widetilde{\mathbb{y}}$.
We can thus construct
the rational tangle~$\widetilde{\mathbb{y}}_r$ associated with the continued
fraction~$r$
with even integers.
The  automorphism~$\mathscr{T}_{\mathbb{y}}$
for  the Dehn twist along~$\mathbb{y}$
is~\eqref{auto_T3}, and~$\mathscr{T}_{\mathbb{x}}$ is given from
the 2-chain relation~\eqref{2-chain} as
\begin{equation}
  \mathscr{T}_{\mathbb{x}}:
  \begin{pmatrix}
    \mathsf{T}_0
    \\
    \mathsf{T}_1
    \\
    \mathsf{X}
    \\
    \mathsf{U}_0
    \\
    \mathsf{U}_1
  \end{pmatrix}
  \mapsto
  \begin{pmatrix}
    \mathsf{X}  \mathsf{T}_0 \mathsf{X}^{-1}
    \\
    \mathsf{T}_1
    \\
    \mathsf{X}
    \\
    \mathsf{X}\mathsf{U}_0\mathsf{X}^{-1}
    \\
    \mathsf{U}_1
  \end{pmatrix} ,
\end{equation}
which is consistent with
the Dehn twist $\sigma_R^2$ on $\Sigma_{0,4}$~\eqref{sigma_R_sphere}.

\begin{figure}[htbp]
  \centering
  $\vcenter{\hbox{\includegraphics[scale=0.7]{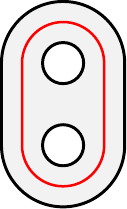}}}
  \overset{\mathscr{T}_{\mathbb{x}}^{-1}}{\longmapsto}
  \vcenter{\hbox{\includegraphics[scale=0.7]{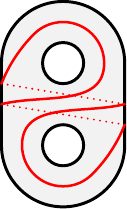}}}
  \overset{\mathscr{T}_{\mathbb{y}}}{\longmapsto}
  \vcenter{\hbox{\includegraphics[scale=0.7]{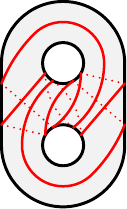}}}
  $
  \caption{A rational tangle for $2+\frac{1}{2}$ which corresponds to
    the figure-eight knot $4_1$.}
  \label{fig:figure-eight}
\end{figure}
We show a few examples.
The figure-eight knot $4_1$ is a rational tangle with
$\frac{5}{2}=2+\frac{1}{2}$ as in Fig.~\ref{fig:figure-eight}.
We have
\begin{equation}
  \label{eq:2}
  \widetilde{\mathbb{y}}_{5/2}
  =
  \left( \mathscr{D}_{\mathbb{y}}
    \mathscr{D}_{\mathbb{x}}^{-1}\right)
  ( \widetilde{\mathbb{y}}) ,
\end{equation}
which gives
\begin{equation}
  \label{eq:3}
  \mathcal{A}(\widetilde{\mathbb{y}}_{5/2})
  =
  \ch
  \left(
    \mathsf{U}_1 \mathsf{T}_0 \mathsf{T}_1
    \left(
      \mathsf{T}_0^\vee
      \mathsf{T}_1 \mathsf{T}_0
    \right)^{-1}
    \mathsf{T}_0^{-1}
    \mathsf{U}_0 \mathsf{T}_1
    \left(
      \mathsf{T}_0^\vee
      \mathsf{T}_1 \mathsf{T}_0
    \right)
  \right) .
\end{equation}
The knot $5_2$  is $\frac{7}{2}=4+\frac{1}{-2}$, and
\begin{equation}
  \label{eq:7}
  \widetilde{\mathbb{y}}_{7/2}=
  \left(
    \mathscr{D}_{\mathbb{y}}^{-1} 
    \mathscr{D}_{\mathbb{x}}^{-2}
  \right)(\widetilde{\mathbb{y}}) ,
\end{equation}
which gives
\begin{equation}
  \mathcal{A}(\widetilde{\mathbb{y}}_{7/2})
  =
  \ch
  \left(
    \mathsf{U}_1 \mathsf{T}_1^{-1} \mathsf{T}_0^{-1}
    \left(
      \mathsf{T}_1
      \mathsf{T}_0
      \mathsf{T}_1^\vee \mathsf{T}_0^{-1}
    \right)
    \mathsf{T}_1 \mathsf{T}_0
    \mathsf{T}_1^\vee
    \mathsf{T}_1 \mathsf{U}_0
        \left(
      \mathsf{T}_1
      \mathsf{T}_0
      \mathsf{T}_1^\vee \mathsf{T}_0^{-1}
    \right)^{-2}
  \right) .
\end{equation}
In both cases,
the DAHA polynomials
$\mathcal{A}(\mathbb{k}_r)(1)$ are too involved to give here.
Nonetheless
Mathematica shows that the constant terms~$\eth^0$ of
$\mathcal{A}(\widetilde{\mathbb{y}}_{r})$ reduce
to
\begin{align*}
  \left.
    \text{Const}(\mathcal{A}(\widetilde{\mathbb{y}}_{5/2}))(1)
  \right|_{x_0=x_1=-x=q^{\frac{1}{2}}}
  & =
    \frac{q^{-2}-q^{-1}+1-q+q^2}{(1-q)(1-q^2)} ,
  \\
  \left.
    \text{Const}(\mathcal{A}(\widetilde{\mathbb{y}}_{7/2}))(1)
  \right|_{x_0=x_1=-x=q^{\frac{1}{2}}}
  & =
  \frac{q(1-q+2q^2-q^3+q^4-q^5)}{(1-q)(1-q^2)} .
\end{align*}
These computations support
a relationship with the Jones polynomial
as observed in~\cite{KHikami19a}.

\section{Concluding Remarks}

We have proposed a generalization of the
type-$C^\vee C_1$ DAHA at
$\mathbf{t}_\star$ by introducing the  Heegaard dual operators.
We hope to report on their roles  on the
(non-symmetric)
Askey--Wilson polynomials at $\mathbf{t}_\star$,
and also
on the generalization to the higher-rank skein algebras.
It would  be promising to  incorporate results from
the cluster algebra~\cite{KHikami17a,ChekhShapi23a}.

\section*{Acknowledgments}
The author would like to thank Hitoshi Murakami
for communications on
rational tangles.
The work of KH is supported in part by
JSPS KAKENHI Grant Numbers
JP22H01117,
JP20K03601,
JP20K03931.


\end{document}